\documentclass[12pt,draft]{amsart}

\makeatletter

\theoremstyle{plain}    
\newtheorem{thm}{Theorem}[section]
\numberwithin{equation}{section} 
\numberwithin{figure}{section} 
\theoremstyle{plain}    
\newtheorem*{thm*}{Theorem} 
\theoremstyle{plain}    
\newtheorem{cor}[thm]{Corollary} 
\theoremstyle{plain}    
\newtheorem{lem}[thm]{Lemma} 
\theoremstyle{plain}    
\theoremstyle{definition}
\newtheorem{defn}[thm]{Definition}
\theoremstyle{remark}

\theoremstyle{remark}

\theoremstyle{remark}    

\theoremstyle{remark}    

\theoremstyle{definition}  

\theoremstyle{remark}
  \newtheorem*{acknowledgement*}{Acknowledgement} 

\theoremstyle{plain}    

\theoremstyle{plain}    
\theoremstyle{plain}    
\theoremstyle{plain}    
\theoremstyle{definition}

\theoremstyle{remark}

\theoremstyle{remark}    

\theoremstyle{remark}    

\theoremstyle{plain}    

\usepackage{amssymb} \usepackage{amscd}

\makeatother

\begin{document}

\title[Popa Algebras]{Excision and a Theorem of Popa}

\author{Nathanial P. Brown}

\address{Department of Mathematics, Penn State University, State
College, PA 16802}

\email{nbrown@math.psu.edu}

\thanks{Supported by an NSF Postdoctoral
Fellowship.  2000 MSC number: 46L05. }

\begin{abstract}
We give an elementary C$^*$-algebraic proof of a result of Sorin Popa which 
is of fundamental importance to Elliott's classification program.  
\end{abstract}

\maketitle

\section{Introduction}

Throughout most of the 1990's, work on the finite case of Elliott's
Classification Program focused on inductive limits of homogeneous
C$^*$-algebras.  Gradually it was noticed that such inductive limits
enjoy a certain internal, finite dimensional approximation property
that is a useful tool in classification (see, for example,
\cite[Theorem 2.21, Corollary 2.24]{EG}).  In 1997 Sorin Popa noticed
that his previous work on injective von Neumann algebras could be used
to show that certain quasidiagonal C$^*$-algebras (cf.\
\cite{brown:QDsurvey}, \cite{dvv:QDsurvey}) enjoy a very similar
internal, finite dimensional approximation property.  These
developments set the stage for Huaxin Lin's definition and subsequent
classification of `tracially AF' algebras (cf.\
\cite{lin:TAFclassification}).  This work is so exciting because
Popa's result only assumes quasidiagonality and deduces an
approximation property which looks similar to the approximation
property which is sufficient for classification.  In other words, one
can now imagine that the mild hypothesis of quasidiagonality would
actually imply an AH inductive limit decomposition! (Of course, other
hypotheses such as nuclearity and simplicity will also be needed.)

While Popa's result is certainly known to many experts it seems that
relatively few have worked through the proof.  This may be partly due
to the fact that the main technical ingredient in the proof (cf.\
\cite[Theorem 2.3]{popa:simpleQD}) relies on some non-trivial
II$_1$-factor theory and a bit of free probability.  In this note we
will give a completely C$^*$-algebraic proof of Popa's result.  Our
hope is that the classification community will find this proof easier
to understand.  However, our real goal is simply to expose and
advertise the beauty and simplicity of Popa's ideas.  Indeed, the
language of this note may differ from Popa's, but the main techniques
and estimates involved come directly from Popa's work.

In honor of Popa's contribution we make the following definition. 

\begin{defn} Let $A$ be a simple, unital C$^*$-algebra.  Then $A$ is
called a {\em Popa algebra} if for every finite set ${\mathfrak F}
\subset A$ and $\epsilon > 0$ there exists a non-zero finite dimensional
subalgebra $B \subset A$ with unit $e$ such that (1) $\| [e,x] \| <
\epsilon$ for all $x \in {\mathfrak F}$ and (2) $e {\mathfrak F} e
\subset^{\epsilon} B$ (i.e.\ for every $x \in {\mathfrak F}$ there
exists $b \in B$ with $\| exe - b \| < \epsilon$).
\end{defn}

The main result in \cite{popa:simpleQD} can then be stated as follows. 

\begin{thm} Let $A$ be simple, unital, quasidiagonal  
(cf.\ \cite{brown:QDsurvey}, \cite{dvv:QDsurvey}) and have real rank
zero (cf.\ \cite{brown-pedersen}).  Then $A$ is a Popa algebra.
\end{thm}

Actually, Popa gets away with less than real rank zero however in the
classification program this is the main case of interest and hence we
will only prove the result in this case.

Before getting into any details we wish to point out that the key
technical results required in this note can be traced back to the
fundamental work of Glimm.  Though these things are now well known, we
state below Glimm's work which we will need. 

\begin{lem}\cite[Lemma 1.10]{glimm:UHF} For every $\epsilon > 0$ and
$n \in {\mathbb N}$ there exists $\delta = \delta(\epsilon,n) > 0$
such that if $A$ is a C$^*$-algebra, $p \in A$ is a projection and $\{
w_1, \ldots, w_n \} \subset A$ is a set such that $\| w_j^* w_i -
\delta_{i,j} p \| < \delta$ ($\delta_{i,j}$ is the Kronecker delta
function) then there exist actual partial isometries $\{ v_1, \ldots,
v_n \} \subset A$ such that $v_j^* v_i = \delta_{i,j} p$ and,
moreover, $\| w_i - v_i \| < \epsilon$ for $1 \leq i \leq n$. 
\end{lem}

This result simply states that a finite set of elements which are
almost partial isometries with common support and orthogonal ranges
can be perturbed to actual partial isometries with common support and
orthogonal ranges.  This is not exactly the content of \cite[Lemma
1.10]{glimm:UHF}, but Glimm's proof is easily adapted to prove the
statement above. 

The other result of Glimm which we will need comes from his 
work on the noncommutative Stone-Weierstrass problem.

\begin{thm}\cite[Theorem 2]{glimm:Stone-Weierstrass} Let $A \subset
B(H)$ be a C$^*$-algebra which acts irreducibly (i.e.\
$A^{\prime\prime} = B(H)$) and which contains no non-zero compact
operators.  Then the pure states on $A$ are weak-$*$ dense in the
state space of $A$.
\end{thm}

\section{Excision of States} 

The main technical idea in the present paper is {\em excision of states}.  
This idea is already present in Popa's paper as \cite[Theorem
2.3]{popa:simpleQD} is really a result about excising factorial
traces.  However, in this section we will observe that general facts
about excision easily give us the tools we need to prove Popa's
result. 

\begin{defn} Let $A$ be a C$^*$-algebra and $\phi \in S(A)$ be a state
on $A$.  We say that {\em $\phi$ can be excised} if there exists a net
of positive, norm one elements $h_{\lambda} \in A$ such that $\|
h_{\lambda}^{1/2} a h_{\lambda}^{1/2} - \phi(a) h_{\lambda}\| \to 0$
for all $a \in A$.  If the $h_{\lambda}$'s can be taken to be
projections then we say {\em $\phi$ can be excised by projections}.
\end{defn}
 
The definitive result concerning excision is due to Akemann, Anderson
and Pedersen: A state can be excised if and only if it belongs to the
weak-$*$ closure of the pure states (cf.\
\cite{akemann-anderson-pedersen}).  When combined with Glimm's results
on density of pure states we arrive at the following result (which is
far from the most general but sufficient for our purposes).

\begin{thm} If $A$ is simple, unital and infinite dimensional then
every state on $A$ can be excised.
\end{thm}

\begin{proof} Since $A$ is simple, unital and infinite dimensional
there exists a (necessarily faithful) irreducible representation whose
image (necessarily) contains no non-zero compact operators.  Hence the
pure states on $A$ are dense in the state space of $A$.  Hence every
state can be excised.
\end{proof}

When there are enough projections around we can even excise by projections. 

\begin{lem} Assume $A$ has real rank zero and $\phi \in S(A)$ can be
excised.  Then $\phi$ can be excised by projections. 
\end{lem}

\begin{proof} Let $h_{\lambda}$ be a net of positive, norm one elements 
in $A$ such that $\| h_{\lambda}^{1/2} a h_{\lambda}^{1/2} - \phi(a)
h_{\lambda}\| \to 0$ for all $a \in A$.  Since $A$ has real rank zero
we may assume that each $h_{\lambda}$ has finite spectrum and hence we
can write
$$h_{\lambda} = \sum_{i = 1}^{k(\lambda)} \alpha_i^{(\lambda)}
Q_i^{(\lambda)},$$ where $\{ Q_i^{(\lambda)} \}$ are orthogonal
projections and $1 = \alpha_1^{(\lambda)} > \alpha_2^{(\lambda)} >
\cdots > \alpha_{k(\lambda)}^{(\lambda)} > 0$.  We claim that the
projections $Q_1^{(\lambda)}$ will excise $\phi$.  Indeed,
we have the following inequality.
\begin{eqnarray*}
\| Q_1^{(\lambda)} a Q_1^{(\lambda)} - 
\phi(a)Q_1^{(\lambda)} \| & = &
\| Q_1^{(\lambda)}\bigg( h_{\lambda}^{1/2} a h_{\lambda}^{1/2} - \phi(a)
h_{\lambda} \bigg) Q_1^{(\lambda)} \|\\
& \leq & \| h_{\lambda}^{1/2} a h_{\lambda}^{1/2} - \phi(a)
h_{\lambda} \|.
\end{eqnarray*}
\end{proof}

\begin{cor} 
\label{thm:excision}
If $A$ is simple, unital, infinite dimensional and has
real rank zero then every state on $A$ can be excised by projections. 
\end{cor}

\section{Popa's Local Quantization Technique}

In \cite{popa:injectiveimplieshyperfiniteI} Popa gave a new proof of
Connes' uniqueness theorem for the injective II$_1$-factor.  The main
idea in the proof is what Popa calls `local quantization' and this
technique has since been exploited with great success.  In this
section we explain how Popa's technique can be used to `excise' certain
finite dimensional, unital, completely positive maps.

The main idea is as follows.  We begin with an algebra $A$ and state
$\phi \in S(A)$ which can be excised by projections.  Let $L^2
(A,\phi)$ be the GNS Hilbert space and $P \in B(L^2 (A,\phi))$ be a
finite rank projection.  Then $PB(L^2 (A,\phi))P$ is just a finite
dimensional matrix algebra and we define a unital, completely positive
map $\Phi: A \to PB(L^2 (A,\phi))P$ by $a \mapsto P\pi_{\phi}(a) P$,
where $\pi_{\phi}$ is the GNS representation.  Popa's local
quantization technique essentially says that in this situation $\Phi$
can be excised.  Note that if a state $\phi$ can be excised by
projections $\{ p_{\lambda} \}$ then we can formulate excision as
follows: There exists a net of $*$-monomorphisms $\rho_{\lambda} :
{\mathbb C} \hookrightarrow A$ (i.e.\ $\alpha \mapsto \alpha
p_{\lambda}$) such that $$\| \rho_{\lambda}(1) a \rho_{\lambda}(1) -
\rho_{\lambda}(\phi(a)) \| \to 0,$$ for all $a \in A$.

If $A$ is a C$^*$-algebra and $\phi \in S(A)$ is a state on $A$ then
for each element $a \in A$ we will let $\hat{a} \in L^2(A,\phi)$
denote the corresponding vector.  

\begin{thm}[Popa's Local Quantization] Let $\phi \in S(A)$ be a state 
which can be excised by projections.  Let $\{ y_i \}_{i = 1}^{m}
\subset A$ be such that $\phi(y_j^* y_i) = \delta_{i,j}$ (i.e.\ in
$L^2(A,\phi)$, $\{ \hat{y}_i \}_{i = 1}^{m}$ is an orthonormal set of
vectors) and $P \in B(L^2(A,\phi))$ be the orthogonal projection onto
the span of $\{ \hat{y}_i \}_{i = 1}^{m}$.  Then there exists a net of
$*$-monomorphisms $\rho_{\lambda} : P\pi_{\phi}(a) P \hookrightarrow
A$ such that $$\| \rho_{\lambda}(P) a \rho_{\lambda}(P) -
\rho_{\lambda}(\Phi(a)) \| \to 0,$$ for all $a \in A$.

Moreover, for each $\lambda$ we have the following commutator
inequality for every unitary element $u \in A$: $$\| [u,
\rho_{\lambda}(P)] \|^2 \leq \| [P, \pi_{\phi}(u)] \|^2 + 2\|
\rho_{\lambda}(P)u\rho_{\lambda}(P) - \rho_{\lambda}(\Phi(u)) \|.$$
\end{thm}

\begin{proof} Let ${\mathfrak F} \subset A$ be a finite
set of norm one elements which, for convenience, contains the unit of
$A$ and $\epsilon > 0$ be arbitrary.  Evidently it suffices to show
that there exists a $*$-monomorphism $\rho : P\pi_{\phi}(a) P
\hookrightarrow A$ such that $\| \rho(P) a \rho(P) - \rho(\Phi(a)) \|
< \epsilon,$ for all $a \in {\mathfrak F}$.

Since $\phi$ can be excised by projections, we can, for any $\delta >
0$, find a projection $p \in A$ such that $\| p(y^*_j x y_i)p -
\phi(y^*_j x y_i)p \| < \delta$ for all $x \in {\mathfrak F}$ and $1
\leq i,j \leq m$.  In particular, note that $\| (y_jp)^* (y_ip) -
\delta_{i,j}p \| < \delta$ for all $1 \leq i,j \leq m$.  In other
words, $\{ y_ip \}_{i = 1}^{m}$ is {\em almost} a set of partial
isometries with orthogonal ranges and common support $p$.

Thus we can perturb the $y_i p$'s to honest partial
isometries $\{ v_i \}$ such that $ v_j^* v_i = \delta_{i,j} p$ and $\|
v_i - y_ip \| < \epsilon/4m^2$ (hence $\delta$ above is a number which
depends on $m$ and $\epsilon$, however, we also assume that $\delta <
\epsilon/2m^2$). Hence if we define $f_{i,j} = v_i v_j^*$ then $\{
f_{i,j} \}$ is a set of matrix units for a $m\times m$-matrix algebra.
Moreover, notice that if we let $q = \sum f_{i,i}$ be the unit of this
matrix algebra then cutting an element $x \in {\mathfrak F}$ we have 
\begin{eqnarray*}
qxq & = & \bigg( \sum_{i = 1}^m f_{i,i} \bigg) x 
          \bigg( \sum_{j = 1}^m f_{j,j} \bigg)\\ 
    & = & \sum_{i,j = 1}^m v_i v_i^* x v_j v_j^*\\
& \approx & \sum_{i,j = 1}^m v_i (y_i p)^* x (y_j p) v_j^*\\ 
& \approx & \sum_{i,j = 1}^m v_i( \phi(y_i^* x y_j)p )v_j^*\\ 
& = & \sum_{i,j = 1}^m \phi(y_i^* x y_j) f_{i,j}.
\end{eqnarray*}
In the two approximations above the norm difference is less than
$\epsilon/2$ and hence the triangle inequality implies that $$\| qxq -
\sum_{i,j = 1}^m \phi(y_i^* x y_j) f_{i,j} \| < \epsilon.$$  

The only thing left to notice is that the matrix of $P \pi_{\phi}(x)
P$ w.r.t.\ the othonormal basis $\{ \hat{y}_i \}_{i = 1}^{m}$ is just
$(<\pi_{\phi}(x)\hat{y}_j, \hat{y}_i> )_{i,j} = (\phi(y_i^* x
y_j))_{i,j}$.  Hence we can identify $PB(L^2(A,\phi))P$ with the
matrix algebra C$^*(\{ f_{i,j} \})$ in such a way that $P
\pi_{\phi}(x) P \mapsto \sum_{i,j = 1}^m \phi(y_i^* x y_j) f_{i,j}$,
for all $x \in A$ and the proof of the first part of the theorem is complete.

Finally, if $x$ is a unitary, the commutator estimate goes as follows.
\begin{eqnarray*}
\| [q,x] \|^2 & = & \| qx - qxq + qxq - xq \|^2\\ 
    & = & \| qxq^{\perp} - q^{\perp}xq \|^2 \\
    & = & \max \{ \|q^{\perp}x^* q\|^2, \| q^{\perp}xq \|^2 \}\\ 
    & = & \max \{ \|qxq^{\perp}x^* q\|, \|qx^* q^{\perp}xq \| \}\\
    & = & \max \{ \|q - qxqx^* q\|, \|q - qx^*qxq \| \}\\
 & \leq & \max \{ \|q - \rho(\Phi(x)\Phi(x^*))\|, \|q - 
               \rho(\Phi(x^*)\Phi(x))\|\} + 2\| qxq - \rho(\Phi(x)) \|\\
    & = & \max \{ \|P - P\pi_{\phi}(x^*)P\pi_{\phi}(x)P\|, 
          \|P - P\pi_{\phi}(x)P\pi_{\phi}(x^*)P\|\} 
          + 2\| qxq - \rho(\Phi(x)) \|\\
    & = & \| [P, \pi_{\phi}(x)] \|^2 + 2\| qxq - \rho(\Phi(x)) \|
\end{eqnarray*}
\end{proof}

We stated the theorem above only for projections onto finite
dimensional subspaces which are actually spanned by vectors coming
from $A$.  However, a simple density argument carries this technique
over to arbitrary finite rank projections.

We are now ready for the harvest.

\begin{thm}\cite[Theorem 3.2]{popa:simpleQD} Let $A$ be simple, 
unital, quasidiagonal and have real rank zero.  Then $A$ is a Popa
algebra.
\end{thm}

\begin{proof} Let $\phi$ be any state on $A$.  Since $A$ is simple and
unital we have that $\pi_{\phi}$ is faithful and its image contains no
non-zero compact operators.  Since $A$ is quasidiagonal we can apply
Voiculescu's Theorem to construct a sequence of non-zero finite rank
projections $P_n \in B(L^2(A,\phi))$ such that $\| [P_n,
\pi_{\phi}(a)] \| \to 0$ for all $a \in A$. 

By Corollary \ref{thm:excision}, $\phi$ can be excised by projections
and hence we can apply Popa's local quantization technique (using the
commutator estimates) to conclude that $A$ is a Popa algebra.
\end{proof}

\bibliographystyle{amsplain}

\providecommand{\bysame}{\leavevmode\hbox to3em{\hrulefill}\thinspace}

\end{document}